\documentclass{article}

\usepackage{amssymb,amsmath}
\usepackage[left=2.5cm,right=2.5cm]{geometry}
\usepackage[english]{babel}

\newcommand{\qed}{\bigskip\hfill\(\Box\)}
\newtheorem{theorem}{Theorem}

\newtheorem{e-proposition}[theorem]{Proposition}
\newtheorem{corollary}{Corollary}
\newtheorem{e-definition}[theorem]{Definition\rm}

\newtheorem{example}{\it Example\/}

\newcommand{\p}{\mathbb{P}}

\newcommand{\E}{\mathbb{E}}
\newcommand{\Var}{\mathrm{Var}}
\newcommand{\R}{\mathbb{R}}

\newcommand*{\ind}[1]{\mathbf{1}_{\{#1\}}}

\newcommand{\PP}{\mathbb P}

\renewcommand{\r}{\right}
\def\la{\left\langle}
\def\ra{\right\rangle}
\newcommand{\bcK}{boundedness condition with constant $K$}

\title{Sharp bounds on the rate of convergence of the empirical covariance matrix\thanks{The research was conducted while the authors participated in the Thematic Program on Asymptotic Geometric Analysis at the Fields Institute in Toronto in Fall 2010.}}

\author{Rados{\l}aw ADAMCZAK\thanks{Research partially supported by MNiSW Grant no. N N201 397437 and the Foundation for Polish Science.} \and Alexander E. LITVAK \and Alain PAJOR \and Nicole TOMCZAK-JAEGERMANN\thanks{This author holds the Canada Research Chair in
  Geometric Analysis.}}

\begin{document}
\maketitle
\begin{abstract}

Let $X_1,\dots, X_N\in\R^n$,
$n\leq N$, be independent centered random vectors
with log-concave distribution and with the identity as covariance
matrix. We show that with overwhelming probability
one has
$$
 \sup_{x\in S^{n-1}} \Big|\frac{1}{N}\sum_{i=1}^N \left(|\langle
   X_i, x\rangle|^2 - \E|\langle X_i, x\rangle|^2\r)\Big|
  \leq C \sqrt{ \frac{n}{N} } ,
$$
where $C$ is an absolute positive constant. This result is valid in
a more general framework when the linear forms $(\langle
X_i,x\rangle)_{i\leq N, x\in S^{n-1}}$  and the Euclidean norms
$(|X_i|/\sqrt n)_{i\leq N}$ exhibit uniformly a sub-exponential
decay. As a consequence, if $A$ denotes the  random matrix with
columns $(X_i)$, then with overwhelming probability, the extremal
singular values $\lambda_{\rm min}$ and $\lambda_{\rm max}$ of
$AA^\top$ satisfy the inequalities $ 1 - C\sqrt{\frac{n}{N}} \le
\frac{\lambda_{\rm min}}{N} \le \frac{\lambda_{\rm max}}{N} \le 1 +
C\sqrt{\frac{n}{N}} $ which is a quantitative version of Bai-Yin
theorem \cite{BY} known for random matrices with i.i.d. entries.
\end{abstract}


\bigskip

Let $X \in \R^n$ be a centered random vector whose
covariance matrix is the identity and $X_1$, ..., $X_N$
be independent copies of $X$.
Let $A$ be a random $n\times N$
matrix whose columns are $(X_i)$.  By $\lambda_{\rm min}$
(resp. $\lambda_{\rm max}$) we denote the smallest (resp. the
largest) singular number of the matrix of empirical covariance
$AA^\top$. In the study of the local regime in the random matrix
theory of particular interest is the limit behavior of extremal
values of the spectrum of $AA^\top$.  In the case of Wishart Ensemble
when the coordinates of $X$ are independent, the Bai-Yin theorem \cite{BY}
establishes the convergence of $\lambda_{\rm min}/N$ and $\lambda_{\rm
  max}/N$ when $n,N\to\infty$ and $n/N\to \beta\in(0,1)$, under the
assumption of a finite fourth moment. In this note we study the
asymptotic non-limit behavior (also called ``non-asymptotic'' in
Statistics) i.e. we look for sharp upper and lower bounds for
singular values in terms of $n$ and $N$, when $n\leq N$ are
sufficiently large. For example, for Gaussian matrices it
is known that singular values satisfy inequalities
\begin{equation}
 1 -C\sqrt{\frac{n}{N}} \le
\frac{\lambda_{\rm min}}{N} \le \frac{\lambda_{\rm max}}{N} \le 1 +
C\sqrt{\frac{n}{N}}
\end{equation}
with probability close to 1. We obtain the same estimates
  for large class of random matrices, which in particular do
not require that entries of the matrix are independent or that
$X_i$'s are identically distributed. Note that the natural question
about convergence of singular values in such a case is still open (see \cite{Au} for the case of $X_i$
having uniform distribution on a rescaled $\ell_p^n$ ball).

\smallskip

The natural scalar product and Euclidean
norm on $\R^n$ are denoted by $\langle\,\cdot, \,\cdot\rangle$ and
$|\cdot|$. We also denote by the same notation $|\cdot|$ the
cardinality of a set. By $C$, $C_1$, $c$ etc. we will denote
absolute positive constants.
%
%
%

Let   $X_1,\ldots,X_N$  be a sequence of random vectors in $\R^n$
(not necessarily identically distributed).
We say that  it is  uniformly $\psi_1$
if for some $\psi > 0$,
\begin{equation}\label{psi_1_assumption}
   \sup _{i\leq N} \sup_{y\in S^{n-1}}
\| \ |\la X_i , y \ra|\ \|_{\psi _1}  \leq \psi,
\end{equation}
where for a random variable $Y\in\R$,
$
\|Y\|_{\psi_1}=\inf\left\{C>0\,;\, \,\E\exp\left({|Y|/C}\right)
\leq 2\right\}.
$
%
%
%
%
%
%
%
We say that
it  satisfies the  \bcK\/ (for some   $K\geq 1$) if
\begin{equation}
  \label{eq:bcK}
  \PP \left( \max _{i\leq N} |X_i|/\sqrt{n} > K  \max\{1,
      (N/n)^{1/4}
\} \right) \leq \exp\left( -  \sqrt{n} \right) .
\end{equation}

\smallskip

The main result of this note is the following theorem.

\begin{theorem} \label{raz}
  Let $n\leq N$ be positive integers and $\psi,
  K \geq 1$. Let $X_1,\ldots,X_N$ be independent random vectors in
  $\R^n$ satisfying (\ref{psi_1_assumption}) and (\ref{eq:bcK}). \/
 Then with probability at least $1 - 2 \exp\left(- c \sqrt{n}\r)$
 one has
$$
 \sup_{x\in S^{n-1}} \Big|\frac{1}{N}\sum_{i=1}^N \left(|\langle
   X_i, x\rangle|^2 - \E|\langle X_i, x\rangle|^2\r)\Big|
  \leq C \left(\psi + K\r)^2 \  \sqrt{ \frac{n}{N} } .
$$
\end{theorem}

\noindent
{\bf Remarks.\ \  1. }
Theorem \ref{raz} improves estimates obtained in
\cite{ALPT} for log-concave isotropic vectors.
There, we considered essentially the case of $N$ proportional to
$n$, which was sufficient to answer the question of Kannan,
Lov\'{a}sz and Simonovits \cite{KLS1}, however, for bigger $N$,
the results were off by a logarithmic factor.
The theorem above removes this factor completely leading to the best
possible estimate for an arbitrary $N$, that is to an estimate of the
same order as in the Gaussian case.
\newline
{\bf 2. }
In the case $N<n$ Theorem~\ref{old} below together with assumptions  (\ref{psi_1_assumption}) and (\ref{eq:bcK}) immediately implies
that the norm of the matrix $A$ with columns $X_1,\ldots,X_N$ satisfies
$$ \|A\| \leq    C \left(\psi + K\r) \sqrt{n}$$
with probability at least $1 - 2 \exp\left(- c \sqrt{n}\r)$.
This in turn implies that $\lambda_{\rm max} \le  C \left(\psi + K\r)^2  n$ and that
$$ \sup_{x\in S^{n-1}} \Big|\frac{1}{N}\sum_{i=1}^N \left(|\langle   X_i, x\rangle|^2 - \E|\langle X_i, x\rangle|^2\r)\Big|
  \leq C \left(\psi + K\r)^2 \  \frac{n}{N}  $$
with the same probability.

\medskip

As a consequence of Theorem~\ref{raz} we obtain  the following quantitative
version of Bai-Yin theorem \cite{BY} known for random matrices with
i.i.d. entries.

\begin{corollary} \label{BaiYin}
Let $A$ be a random $n \times N$ matrix, whose columns $X_1,\ldots,X_N$
are isotropic random vectors satisfying the assumptions of Theorem \ref{raz}.
%
%
Then with probability at least $1 - 2\exp(-c\sqrt{n})$,
$$
1 -  C \left(\psi + K\r)^2\sqrt{\frac{n}{N}} \le
\frac{\lambda_{\rm min}}{N} \le \frac{\lambda_{\rm max}}{N} \le 1 +  C
\left(\psi + K\r)^2\sqrt{\frac{n}{N}}.
$$
\end{corollary}

\medskip

To emphasize the strength of the above results we observe that
conditions (\ref{psi_1_assumption}) and (\ref{eq:bcK}) are valid for
many classes of distributions.


\begin{example}
  \label{example-one} {\rm Random vectors uniformly distributed on the
    Euclidean ball of radius $K \sqrt{n}$ clearly satisfy (\ref{eq:bcK}).
They also satisfy (\ref{psi_1_assumption}) with $\psi = CK$.}
\end{example}

\begin{example}
  \label{example-two}
{\rm Log-concave isotropic random vectors in $\R^n$.
Recall that a random vector is isotropic if its covariance matrix is
the identity and it is log-concave if its distribution has a
log-concave density.
Such vectors satisfy (\ref{psi_1_assumption}) and (\ref{eq:bcK}) for
appropriate absolute constants $\psi$ and $K$. The boundedness
condition follows from Paouris' theorem (\cite{Pao}) and is
explicitly written e.g., in \cite{ALPT}, Lemma 3.1. We would like to remark that a version of Theorem
\ref{raz} with a weaker probability estimate was proved by Aubrun in the case of isotropic
log-concave random vectors under an additional assumption of
unconditionality (see \cite{Au1}).
%

}
\end{example}

\begin{example} {\rm Any isotropic random vectors $(X_i)_{i\le N}$ in $\R^n$,
satisfying the Poincar\'{e} inequality with constant $L$, i.e. such
that $\Var(f(X_i)) \le L^2 \E |\nabla f(X_i))|^2$ for all compactly
supported smooth functions, satisfy (\ref{psi_1_assumption}) with
$\psi = CL$ and (\ref{eq:bcK}) with $K = CL$.
The question from \cite{KLS} whether all log-concave isotropic random vectors
satisfy the Poincar\'{e} inequality with an absolute constant is one of the major open
problems in the theory of log-concave measures}.
\end{example}

The proof of Theorem \ref{raz} is close to  arguments in Section 4.3 of
\cite{ALPT}, however
it uses a choice of parameters more appropriate for the case
considered here,  and a new approximation argument.
%
We need additional notations.  Let $1\le m\leq N$. By $U_m$ we denote
the subset of all vectors in $S^{N-1}$ having at most $m$ non-zero
coordinates. For an $n \times N$ matrix $A$ we let
 \begin{equation}
  \label{eq:A_m}
  A_m =  \sup _{z\in {U_m}} |Az|.
\end{equation}
The main technical tool is the following result which is the ``in
particular" part of Theorem~3.13 from \cite{ALPT} in which one needs
to adjust corresponding constants and to take a union bound.

\begin{theorem}
\label{old}
Let $X_1,\ldots,X_N$ be as in Theorem
  \ref{raz}, let $A$ be a random $n\times N$ matrix whose columns are
  the $X_i$'s.
Then for every $t\geq 1$ one has
$$
  \PP \left( \exists m \ \  A_m \geq C \psi t
              \max\{\sqrt{m} \ln \frac{2N}{m}, \sqrt{n}\}
  + 6 \max _{i\leq N} |X_i|  \r)
   \leq  \exp\left( -  t \sqrt{n} \r).
$$
\end{theorem}







\noindent
{\bf Proof of Theorem~\ref{raz}. }
%
%
  For $x \in S^{n-1}$ set
$$
  S(x) = \left|\frac{1}{N} \sum_{i=1}^N \left( |\langle
  X_i, x \rangle|^2  - \E |\langle X_i, x\rangle|^2  \r) \r|.
$$
Let $B>0$ be a parameter which we specify later and observe that
\begin{eqnarray*}
\lefteqn{
   \sup_{x\in S^{n-1}} S(x) \le  \sup_{x\in S^{n-1}} \Big( \Big|\frac{1}{N}
  \sum_{i=1}^N \left( \left(|\langle X_i, x\rangle|\wedge B \r)^2  -
  \E \left(|\langle X_i, x\rangle|\wedge B \r)^2  \r) \Big| \Big.
}\\
&+& \Big.\frac{1}{N} \sum_{i=1}^N \left( |\langle
  X_i, x\rangle|^2  - B^2 \r) \ind{|\langle X_i, x\rangle|\ge B} +
   \frac{1}{N}\E \sum_{i=1}^N \left( |\langle
  X_i, x\rangle|^2  - B^2 \r) \ind{|\langle X_i, x\rangle|\ge B} \Big).
\end{eqnarray*}
We denote the summands under the supremum by $S_1(x)$,
$S_2(x)$, and $S_3(x)$, respectively.


\medskip

\noindent
{\bf Estimate for $S_1$:\ }
 Given
$x\in S^{n-1}$ and $i\leq N$ let $Z_i = Z_i (x) = \left(|\langle X_i, x\rangle|\wedge B
\r)^2 - \E \left(|\langle X_i, x\rangle|\wedge B \r)^2$.
Then $|Z_i| \leq B^2$. Moreover,  since
$$
  \Var(Z_i)\leq \E\left( |\langle X_i, x\rangle| \wedge B\r)^4
  \leq \E |\langle X_i, x\rangle| ^4  \leq C_1 \psi ^4,
$$
we observe that
$
 \sigma^2 = \frac{1}{N} \sum_{i=1}^N \Var(Z_i)\leq C_1 \psi ^4
$.
Thus, by Bernstein's inequality
$$
 \PP \left( S_1(x) \geq \theta  \r) =
 \PP \left( \frac{1}{N} \sum_{i=1}^N Z_i \geq \theta  \r) \leq
 \exp\left(-\frac{\theta^2 N}{2( C_1 \psi ^4  + B^2 \theta/3) }\r).
$$


\noindent  It is well known that $S^{n-1}$ admits a $(1/3)$-net
$\cal{N}$ in the Euclidean metric such that $|{\cal N}| \leq 7^n$.
Then by the union bound we obtain that if
\begin{equation}\label{cond3}
  \theta^2 N > 8 C_1 \psi ^4 n \ln 7
 \quad \mbox{ and } \quad
  \theta N > (8/3) B^2 n \ln 7
\end{equation}
then
\begin{equation}\label{s-one}
  \PP\left( \sup _{x\in {\cal{N}}} S_1 (x) \geq \theta \r) \leq
 \exp\left(n\ln 7 -\frac{\theta^2 N}{2( C_1 \psi ^4  + B^2 \theta/3) }\r)\le
  \exp\left( - \frac{\theta^2 N}{4( C_1 \psi ^4 + B^2 \theta/3) }\r).
\end{equation}




\medskip

\noindent
{\bf Estimates for $S_2$ and $S_3$:\ } By H\"older's inequality and
(\ref{psi_1_assumption}) we have, for some absolute constant $C_2 \ge 1$,
\begin{equation}\label{S3}
\sup_{x\in S^{n-1}} S_3(x) \le \frac{1}{N}\sum_{i=1}^N \sup_{x\in S^{n-1}}\|\langle X_i,
x\rangle\|_4^2\ \p\left( |\langle X_i, x\rangle| \ge B\r)^{1/2} \le C_2 \psi^2\exp(-B/\psi).
\end{equation}

To estimate $S_2$, we will use the following notation
$$
  M= \max\bigl\{ \psi^2 n, \max _{i\leq N} |X_i|^2\bigr\}, \quad
 E_B = E_B(x) = \{ i\le N \colon |\langle X_i, x\rangle|\ge B\},\quad
m= \sup _{x\in S^{n-1}} |E_B (x)| .
$$
%
By the definition of $A_{m}$,  we have for every $x\in S^{n-1}$
$$
  B^2 |E_B| \leq \sum_{i \in E_B}|\langle X_i, x\rangle|^2
  \le  \sup _{|E| \leq m} \sum _{i\in E}
  |\langle X_i, x\rangle|^2 \leq A_{m}^2 ,
$$
which yields
$  B^2  m\leq {A^2_{m}} $ and $NS_2(x)\leq A_m^2$.
Theorem~\ref{old}
implies that
for some absolute constant $C\geq C_2$,
with probability at least $1-\exp(-\sqrt{n})$ one has
\begin{equation}\label{B-est}
  B^2 m\leq C \left(M + \psi^2 m \ln^2 \frac{2 N}{m}\r)
  \quad {\rm and}\quad
  \sup _{x\in S^{n-1}}S_2(x) \leq C \left(\frac{M}{N} + \psi^2 \frac{m}{N}
  \ln^2 \frac{2 N}{m}\r).
\end{equation}
Now we choose the parameters. Let $B= 2\sqrt{2C} \psi \ln (5 N/n)$.
Then (\ref{S3}) gives
$
S_3(x) \le C\psi^2 \frac{n}{N} \le
C\frac{M}{N}$ for all $x\in S^{n-1}$ and together with (\ref{B-est}) it yields
that with probability at least $1-\exp(-\sqrt{n})$ one has
$$
  \sup _{x\in S^{n-1}} (S_2(x) + S_3(x)) \leq C \left((2{M}/{N}) +
  \psi^2({m}/{N})\ln^2 ({2 N}/{m})\r) .
$$

It is easy to check that $M \geq \psi^2m \ln^2 ({2 N}/{m})$ on the
set where (\ref{B-est}) holds. Indeed,
assume it is not so, thus $  M<\psi^2 m \ln^2 ({2N}/{m})$.
Then by (\ref{B-est}) we observe  that $B^2 \leq 2 C \psi^2 \ln^2
({2 N}/{m})$, which implies
$$
   m\leq 2 N \exp(-B/\psi \sqrt{2C}) = {2 n^2}/{25 N}.
$$
By our hypothetical upper bound for $M$  and since $f(m)=m\ln^2(2N/m)$
increases on $[1, 2N/e^2]$, we get
$$
 \psi^2 n\leq M\leq \psi^2 ({8 n^2}/{25 N}) \ln^2 ({5 N}/{n}),
$$
which is impossible.


It follows that
$$
  \PP \left( \sup _{x\in S^{n-1}} (S_2(x) + S_3(x)) \leq 3 C ({M}/{N}) \r)
  \geq 1- \exp(-\sqrt{n}).
$$
Combining this estimate with (\ref{s-one}), we get
$$
 \PP \left(\sup_{x\in {\cal{N}}} S(x) \leq \theta + 3 C\frac{M}{N}\r)
 \geq  1- \exp(-\sqrt{n}) -
  \exp\left( - \frac{\theta^2 N}{4(C_1 \psi ^4 + B^2 \theta/3) }\r).
$$
We now  set $  \theta  = C_3  \psi^2 \ \sqrt{n/N}$,
where
$C_3$ is a sufficiently large absolute positive constant so that
(\ref{cond3}) is satisfied. Then using \bcK\/  we obtain
$$
 \PP \left(
 \sup_{x\in {\cal{N}}} S(x) \leq \left( C_3  \psi^2 + 3 C K^2 \r)
 \sqrt{n/N}
 \r) \geq  1- \exp(-\sqrt{n}) - \exp\left( - c n \r)
  \geq 1- 2\exp(- c \sqrt{n}),
$$
where $c$ is a sufficiently small positive constant.
It proves the desired estimate on the $(1/3)$-net.


To pass from ${\cal{N}}$ to the whole sphere  note that $S(x)$
can be written as $|\langle T x, x\rangle |$, where $T$ is a
self-adjoint operator on $\R^n$.  Thus, writing for each $x \in
S^{n-1}$, $x = y +z$ with $y \in {\cal N}$ and $|z|\le 1/3$, we get
\begin{displaymath}
  \|T\| = \sup_{x\in S^{n-1}} |\langle T x, x\rangle |
\le
\sup_{y\in {\cal{N}}} |\langle Ty, y \rangle | +
\frac{2}{3} \sup_{y\in {\cal{N}}} | Ty | +
\sup_{|z|\le 1/3} |\langle Tz, z \rangle |
   \leq
   \sup_{y\in {\cal{N}}} S(y)+ \frac{7}{9} \|T\| ,
\end{displaymath}
which implies the desired  estimate on the whole sphere $S^{n-1}$.
\qed


\begin{flushleft}
Rados{\l}aw Adamczak\\
Institute of Mathematics,\\
University of Warsaw\\
Banacha 2, 02-097 Warszawa, Poland \\
{\tt R.Adamczak@mimuw.edu.pl}

\smallskip
Alexander E. Litvak\\
Department of Mathematical and Statistical Sciences,\\
University of Alberta,\\
Edmonton, Alberta, Canada T6G 2G1\\
{\tt alexandr@math.ualberta.ca}

\smallskip
Alain Pajor\\
Equipe d'Analyse et Math\'ematiques Appliqu\'ees,\\
Universit\'e Paris Est, \\
5 boulevard Descartes, Champs sur Marne, 77454 Marne-la-Vallee, \\
Cedex 2, France\\
{\tt alain.pajor@univ-mlv.fr}

\smallskip
Nicole Tomczak-Jaegermann,\\
Department of Mathematical and Statistical Sciences,\\
University of Alberta,\\
Edmonton, Alberta, Canada T6G 2G1\\
{\tt nicole@ellpspace.math.ualberta.ca}
\end{flushleft}
\end{document}